\documentclass{compositio}
\usepackage{epsf,latexsym,graphicx}
\usepackage{psfrag}
\usepackage{amsmath}

\theoremstyle{plain}

\newtheorem{theorem}{Theorem}[section]
\theoremstyle{definition}
\newtheorem{definition}[theorem]{Definition}
\theoremstyle{lemma}
 \newtheorem{lemma}[theorem]{Lemma}
\theoremstyle{corollary}
 \newtheorem{corollary}[theorem]{Corollary}
\theoremstyle{proposition}
\newtheorem{proposition}[theorem]{{Proposition}}
\theoremstyle{clam}
 
\theoremstyle{conjecture}
\newtheorem{con}[theorem]{Conjecture}
\theoremstyle{remark}

\newtheorem{problem}[theorem]{Problem}

\newtheorem{examp}[theorem]{Example}

\newcommand{\tensor}{\otimes}

\newcommand{\Proj}[1]{\mathbb{ P}(#1)}

\newcommand{\real}{\mathbb R}

\newcommand{\zed}{\mathbb Z}

\newcommand{\oo}{\mathcal O}

\newcommand{\ool}[1]{{\mathcal O}_{#1}}
\newcommand{\oof}[2]{{\mathcal O}_{#1}({#2})}

\newcommand{\Pin}[1]{{\mathbb P}^{#1}}

\newcommand{\ex}{X(\Delta)}

\newcommand{\exy}{X(\Delta')}
\newcommand{\exdual}{X(\Delta)^*}
\newcommand{\exd}{X(P)^*}

\newcommand{\el}[1]{{\mathcal L}_{#1}}

 \newcommand{\proj}[1]{\mathbb{P}^{#1}}
\newcommand{\e}{{\bf e}}

\newcommand{\s}{\overline{s}}
\newcommand{\mult}{\mathrm{mult}}

\begin{document}

 \title{Toric manifolds with degenerate dual variety and defect polytopes}

\author[Sandra Di Rocco]{Sandra Di Rocco}
\thanks{The work for this article was carried out   at Yale University and at the University of Minnesota. The author's stay at Yale was partially supported by the G\"oran Gustafsson Stiftelse. Both visits were possible thanks to a leave of absence granted by KTH, Sweden.} 
\address{Department of Mathematics\\Royal Institute of Technology, S-10044 Stockholm, Sweden}
\curraddr{ University of Minnesota. Math Dept., 127 Vincent Hall, 206 Church St. S.E
Minneapolis, MN 55455 USA}

 \begin{abstract}
Toric manifolds  with dual defect  are classified.
The associated polytopes, called defect polytopes, are proved to be
the class of Delzant integral polytopes for which a combinatorial invariant vanishes.

\end{abstract}

\maketitle


\section{Introduction}
In this paper we address the following two apparently disjoint questions, one of a purely
complex geometrical nature and one of a purely combinatorial nature.

Let $\ex\subset\proj{n}$ be a projective toric manifold.
The dual variety $\ex^*$ is the algebraic variety in $\proj{n^*}$ parameterizing the hyperplanes singular at $X$. Since there is a $(codim(\ex)-1)$-dimensional space of hyperplanes tangent to each point of $X$,  the dual variety $\ex^*$ is expected to be a hypersurface. If the dimension is less than expected the toric manifold is said to have  dual defect.

It is well known that  there are no curves with dual defect and that in dimension  two the only projective manifold having dual defect is the projective plane.

\begin{problem}\label{q1} Characterize the projective toric manifolds
  whose dual variety has lower dimension than expected, namely the
  toric projective manifolds with positive {\em dual defect} $d(\ex)=n-1-dim({\ex}^*)$. 
\end{problem}

Using properties of the toric action and a convenient application of
results in \cite{bfs} we give a complete classification of toric
manifolds with dual defect. In \S\ref{class} we prove that an
$r$-dimensional toric manifold  $\ex$ has positive dual defect $d(\ex)>0$ if and only if:
\begin{center}

\begin{tabular}{llll}\label{table}
$\bullet$&$\ex=\proj{r}$& $d(\ex)=r$&  \\
&&&\\ 
$\bullet$& $\ex=
\begin{array}{c}\Proj{\el{0}\oplus...\oplus\el{\frac{r+d}{2}}} \\
\downarrow \\Y(\Delta') \end{array}$& $d(\ex)=\left\{ \begin{array}{cc}
2,4,...,r-2&if \,\,r\,\, even\\
1,3,...,r-2&if \,\,r\,\, odd\end{array}\right. $&\\
&&&
\end{tabular}
\end{center}
where $\el{i}\in Pic(Y(\Delta'))$ and  $Y(\Delta')$ is a  $\frac{r-d}{2}$-dimensional toric manifold.  

\vspace{.1in}

 Let $P$ be a convex integral polytope of dimension $r$. Define the following combinatorial invariant:
\begin{equation}\label{c(P)}
c(P):=\sum_{F\in F(P)}(-1)^{\mathrm{codim}(F)}(\mathrm{dim}(F)+1)!Vol(F)
\end{equation}
where $F(P)$ is the set of nonempty faces of $P$, including $P$ itself, the volume is the integral volume (see \S\ref{defpoly}) and the volume of a vertex is set to be $1$.
For example the polytope $P=\Delta_2\times I$, the product of a two dimensional standard simplex and the unit interval (see Figure \ref{poly}) has 
{$$c(P)=4!\frac{1}{2}-3!4+2\cdot 9-6=0$$
\begin{figure}
\includegraphics{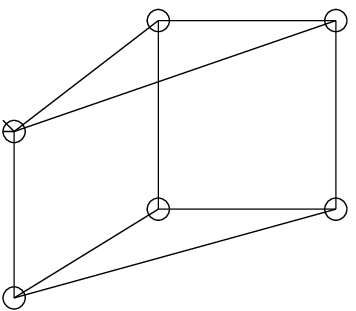}
\begin{caption}{\label{poly}$P=\Delta_2\times I$ in \ref {q2}}
\end{caption}
\end{figure}

\begin{problem}\label{q2} It is natural to ask if $c(P)$ is in fact always a non negative integer and if so try to characterize the polytopes achieving the minimum, i.e. for which $c(P)=0$.
\end{problem}
In \S\ref{nonneg} we prove that for a convex integral simple polytope
a related invariant $c^*(P)$, which coincides  with $c(P)$ if the
toric variety is non singular, is nonnegative.  The nonnegativity of
$c(P)$ for simple polytopes was already established in  \cite{GKZ} using a different approach, see \S\ref{conjsec} for more details.

In \S\ref{defpoly} we characterize the polytopes associated to a non singular toric variety, achieving the minimum. One can find several names for them in the literature. Algebraic geometers refer to them as {\it very simple} polytopes (\cite{Fu}) or {\it absolutely simple } polytopes (\cite{Oda}). Symplectic geometers   use the term {\it Delzant } polytopes. We will use the latter terminology.

In \S\ref{dualpoly} we show that, for a Delzant polytope, $c(P)$ is zero if and only if $P$ is one of the following:
\begin{center}
\begin{tabular}{llll}\label{table}

$\bullet$&$P=\Delta_r$&the standard $r$-dimensional simplex&  \\
&&&\\ 
$\bullet$& $P=\Proj{P_0,...,P_{\frac{r+d}{2}}}$&
 $d=\left\{ \begin{array}{cc}
2,4,...,r-2&if \,\,r\,\, even\\
1,3,...,r-2&if \,\,r\,\, odd\end{array}\right. $&
\end{tabular}
\end{center}
where $P_i$ are Delzant, $\mathrm{dim}(P_i)=\frac{r-d}{2}$ and $\Proj{P_0,...,P_{i}}$ is the projective join defined in \S \ref{join}. 

\vspace{.1in}
 
We show that the  two problems are equivalent under the standard dictionary between toric geometry and convex polytopes. Hence we can use the 
characterization from \ref{q1} to give a complete answer to \ref{q2}.
\acknowledgements{Section \ref{conjsec} contains considerations which were suggested to me by Victor Reiner and Hugh Thomas. I am very grateful to them for their interest in this problem and for sharing their ideas and knowledge on the subject. A special thanks to  Jeremy Martin for several conversations on convex polytopes and to Alicia Dickenstein for pointing out to me the reference \cite{GKZ}.}

\section{Toric manifolds with degenerate dual variety}\label{dual}
Let $\ex\subset\proj{n}$ be a non singular $r$-dimensional projective toric variety (toric manifold) associated to the fan $\Delta$. We will use the geometrical terminology and call such a fan {\it regular}. Combinatorists typically use the term {\it unimodular}.

For basic definitions and properties of toric varieties we refer to \cite{Oda, Fu}. 

Let ${\proj{n}}^*$ be the dual projective space. The associated dual variety is defined as  
$$\exdual=\overline{\{H\in {\proj{n}}^*|H \text{ is singular at some }x\in\ex\}} $$
Generally one expects to impose $(\mathrm{codim}(\ex)-1)$ conditions  for a hyperplane to be tangent at a point on $X$. Hence one  expects the dual variety to be a hypersurface in ${\proj{n}}^*$. 

If the variety is a curve this is indeed the case. For
surfaces the dual variety is a hypersurface unless $X=\proj{2}$. In
higher dimension there are more exceptions. 
 
The dual defect of a variety $X$ is defined as
$d(X)=n-1-dim(X^*)$. In \cite{ls} and \cite{book} nonsingular projective varieties of dimension up to $10$ with positive dual defect  are classified. Fundamental results on dual varieties can be found in \cite{E1,E2}.

 We will call a toric projective manifold a {\em defect toric manifold} if $d(\ex)>0$.

In the following two lemmas we are going to collect some well known results
which will be essential in the proof of the main proposition.
\begin{lemma}\label{proj} Let $\ex$ be an $m$-dimensional  defect toric manifold, with dual defect $d>0$. If the Picard group of $\ex$ has rank $1$, then $\ex\cong \proj{m}$ and $d=m>2$.
\end{lemma}
\begin{proof} The only toric manifold with Picard group of rank one is $\proj{m}$, since it is the only toric manifold represented by a fan consisting of $m+1$ edges. The fact that it has positive defect implies that $\ex$ is a linear $\proj{m}$ and that $d=m$, see \cite[1.6.12]{book}.\end{proof}

Let $\tau_X=min\{t\in\real\,|\, K_{\ex}\otimes\oof{\ex}{t}\text{ is
  nef }\}$.   Because every nef line bundle is  globally generated on
a toric manifold the linear system $|K_{\ex}\otimes\oof{\ex}{\tau}|$
defines a morphism \\
$g:\ex\to\proj{M}$. Let $g=s\circ \psi_\tau$ be the Remmert-Stein factorization, where $\psi_\tau:\ex\to Y$ is a morphism with connected fibers and a normal image and $s$ is a finite-to-one morphism  onto $Im(g)$.  The morphism $\phi_{\tau}$ is the {\it nef value morphism}, defined for example in \cite[1.5]{book}.
\begin{lemma}\label{morph} Let $\ex$ be a defect toric manifold. Then there is a morphism with connected fibers $\psi:\ex\to Y$, where $Y$ and a general fiber $F$ are toric varieties.\end{lemma}
\begin{proof}Consider the nef-value morphism $\psi_{\tau}$, associated to $\ex$. It is a classical fact, \cite{lib}, that a morphism with connected fibers between normal projective varieties, induces a homomorphism from the connected component of the identity of the automorphism group of $\ex$ to the connected component of the identity of the automorphism group of $Y$, with respect to which $\psi_{\tau}$ is equivariant. This implies that both $Y$ and a general fiber are toric.\end{proof}

  These simple observations and the results in \cite[1.2]{bfs}, \cite[14.4]{book} yield the following.

\begin{proposition}\label{dual}
Let $\ex$ be a defect toric manifold of dimension $r$ and defect $d$.  Then there is a non singular toric variety $Y$ of dimension $\frac{r-d}{2}$,  and  invariant line bundles $\el{0},...,\el{\frac{r+d}{2}}$ over $Y$ such that $\ex=\Proj{\el{0}\oplus...\oplus\el{\frac{r+d}{2}}}$. Moreover $r+d\geq 4$.
\end{proposition}

\begin{proof}
Let $\psi:\ex\to Y$ be the morphism described in \ref{morph}.
Since $d(\ex)>0$ \cite[14.4.3]{book} implies that $Pic(F)\cong\zed$
for a general fiber $F$ of the nef-value morphism $\psi_\tau=\psi$. Since
$F$ is toric and non singular, it follows by \ref{proj} that $F\cong \proj{r-m}$. 

Moreover considering $F$ as a projective variety via the embedding given by
$\oof{\ex}{1}|_F$ one has $d(F)>0$, which implies that $\oof{\ex}{1}|_F\cong\oof{\proj{r-m}}{1}$. Then  $\psi_\tau$ is a morphism onto a normal variety such that for a general fiber $F$, $(F,\oof{\ex}{1}|_F)\cong(\proj{r-m},\oof{\proj{r-m}}{1})$. It follows
that, see for example \cite[3.2.1]{book}, $\ex$ has the structure of a linear $\proj{r-m}$-bundle, $\ex\cong\Proj{{\psi_\tau}_*\oof{\ex}{1}}$. Because $\ex$ is non singular, the toric variety $Y$ is also non singular.

Since $\ex$ is toric the vector bundle ${\psi_\tau}_*\oof{\ex}{1}$
splits into the sum of line bundles, see \cite[1.1]{chern}. Hence $\ex=\Proj{\el{0}\oplus...\oplus\el{k}}.$

In \cite{bfs} it is proven that $\tau=\frac{r+d}{2}+1$ and thus $K_F\otimes\oof{F}{\frac{r+d}{2}+1}=\oo$. This implies that $dim(Y)=\frac{r-d}{2}$ and $dim(F)=\frac{r+d}{2}$. 

Finally  since $d(F)>0$, $dim(F)=\frac{r+d}{2}\geq 2$.\end{proof}

It is classically known, see for example \cite[5.2]{ls1}, that a $\proj{k}$-bundle over a variety of dimension $m\leq k$ has dual defect $d=k-m$.

From parity reasons one sees that:

\begin{corollary}\label{class}

A toric manifold $\ex$ of dimension $r$ is a defect toric manifold if and only if it is one of the following:
\begin{itemize}
\item  $\proj{r}$ or $\ex=\Proj{\el{0}\oplus...\oplus\el{2k}}$ over a toric
manifold  $Y$ of dimension $\geq 1$, for $k=1,..., \frac{r-2}{2}$ \\
 ($r$ must be even). The defect is $d(\ex)=r,2k$, respectively.
\item  $\proj{r}$ or $\ex=\Proj{\el{0}\oplus...\oplus\el{2k+1}}$ over a toric
manifold  $Y$ of dimension $\geq 1$, for $k=0,..., \frac{r-3}{2}$ ($r\geq 3$
  must be odd). The defect is $d(\ex)=r,2k+1$, respectively.
\end{itemize} 

\end{corollary}

 In particular we get classically known results like:

\begin{itemize}
\item the only two dimensional defect toric manifold is $\proj{2}$.
\item the only $r$-dimensional toric manifold with maximum defect $d=r$ is $\proj{r}$.
\item the only $r$-dimensional toric manifolds with defect $r-1$ are $\proj{r-1}$-bundles over $\proj{1}$.
\end{itemize}
We will now describe the projective structure of a toric defect $r$-dimensional manifold, namely the invariant subvarieties of an equivariant $\proj{k}$-bundle on a non singular $m$-dimensional toric manifold $Y=\exy$, $\pi:\ex\to\exy$. 
\subsection{The tautological bundle}
The fan $\Delta$ is obtained from $\Delta '$ in the following way, see \cite[1.33]{Oda} for details. 
Let $\{e_1,...,e_k\}$ be a $\zed$-basis for $\real^k$ and let $e_0=-e_1-...-e_k$. Denote by $\sigma^i$ the cone $\oplus_{j\neq i}(\real_+e_j)$. Let $\{n_1,...,n_m\}$ be a $\zed$-basis for $\real^m$. 

Let $Div(X(\Delta'))=\oplus \zed \cdot D_i$, where $D_i$ are the invariant principal divisors associated to the edges $n_i$ in the fan $\Delta'$.
Assume $\ex=\Proj{\el{0}\oplus...\oplus\el{k}}$, where $\el{j}=\sum_ia^j_iD_i$.

Define a map $p:\real^{m}\to \real^m\oplus\real^k$  by $p(n_i)=(n_i, \sum_0^k a^j_ie_j)$. Then $$\Delta=\{p(\sigma)+\sigma^i|\sigma\in\Delta', i=0,...,k\}$$
Let $\Delta(s)$ denote the collection of the $s$-dimensional cones in $\Delta$ and let $$V:\Delta(s)\to\{\text{codimension }s\text{ invariant subvarieties in }X\}$$ denote the $1-1$ correspondence between $s$-dimensional cones of $\Delta$ and $(r-s)$-dimensional invariant subvarieties of $\ex$. We will use the following notation:
\begin{center}
\begin{tabular}{ccc}\label{subvarieties}
$s$&$\sigma\in\Delta(s)$&$V(\sigma)$\\

\hline
\vspace{.05in}
$r$&$p(\sigma)+\sigma^i$, for $\sigma\in\Delta'(m)$& fixed point $m_{\sigma}^i$\\
\vspace{.05in}
$1$&$\e_i=e_i$&  $V(\e_i)$\\
\vspace{.05in}
&$w_i=p(n_i)$&$V(w_i)$\\
\vspace{.05in}
$r-1$&$\rho^i=p(\rho)+\sigma^i$, for $\rho\in\Delta'(m-1)$&$V(\rho^i)$\\
\vspace{.05in}
&$\rho_{ij}=p(\sigma)+\sigma^i\cap\sigma^j$, for every $\sigma\in\Delta'(m)$&$V(\rho_{ij}^{\sigma})$
\end{tabular}
\end{center}
\vspace{.1in}
In order to simplify notation the symbol $n_i$ is used simultaneously for the base vector and the corresponding ray.

 It follows that:
$$Div(\ex)=\left(\bigoplus_{n_i\in\Delta'(1)}\zed\cdot V(w_i)\right)\bigoplus\left(\bigoplus_{j=0,...,k}\zed\cdot V(\e_j)\right)$$
Moreover the linear equivalences among the generators are:
$$-V(\e_0)+V(\e_1)+\sum_{n_i\in\Delta'(1)}(a^i_1-a^0_1)V(w_i)=0$$
 Observing that $\pi^* V(n_i)=V(w_i)$ yields
\begin{equation}\label{equivalence}
V(\e_i)=V(\e_0)+\pi^*(\el{0})-\pi^*(\el{i})
\end{equation}

On the other hand for any projective bundle $\pi:\Proj{E}\to Y$ of rank $k$ one has
$$Pic(\Proj{E})=\zed\cdot\xi_E\oplus\pi^*(Pic(Y))\;\;\;\text{and}$$
\begin{equation}\label{classic}
K_{\Proj{E}}=\pi^*(K_Y+ det(E))-k\xi_E
\end{equation}
where $\xi_E$ is the tautological line bundle associated to $E$. For simplicity of notation let us denote by $\xi$ the tautological bundle associated to $\el{0}\oplus...\oplus\el{k}$. The description of the canonical bundle in terms of the principal invariant divisors gives:
\begin{equation}\label{toric}
K_{\ex}=-\sum_{n_i\in\Delta'(1)} V(w_i)-\sum_{j=0,...,k} V(\e_j)\;\text{ and }\;K_{\exy}=-\sum_{n_i\in\Delta'(1)} V(n_i)\end{equation}
Here we are using an additive structure to follow a standard toric notation. Comparing (\ref{classic}) and (\ref{toric}) we obtain
$$-\sum_{n_i\in\Delta'(1)}V(w_i)+\sum_1^k \pi^*(\el{i})-k\cdot\xi=-\sum_{n_i\in\Delta'(1)} V(w_i)-\sum_{j=0,...,k} V(\e_j)$$
The linear relations (\ref{equivalence}) yield:
$$\sum_1^k \pi^*(\el{i})-k\cdot\xi=-kV(\e_0)-k\pi^*(\el{0})+\sum_{j=1,...,k}\pi^*(\el{i})   $$
which proves that:
\framebox{$\xi=V(\e_0)+\pi^*(\el{0})$}
\subsection{Intersection numbers}\label{polytope}
In order to understand the positivity of a line bundle on a toric variety one has to give an estimate of the intersection numbers of the line bundle with all the invariant rational curves.
In particular let $L=\sum a_iD_i$ be a line bundle on an $r$-dimensional nonsingular toric variety $\ex$ and let $V(\rho)$ be a rational invariant curve corresponding to the $(r-1)$-dimensional cone $\rho=\sigma_1\cap\sigma_2$. Then if $\sigma_1=<v_1,v_2,...,v_r>$ and $\sigma_2=<v_0,v_2,v_3,...,v_r>$ the intersection number $L\cdot V(\rho)$ is given by
$$L\cdot V(\rho)=a_0+a_1-\sum_2^r s_ia_i$$
where the $s_i$'s are the integers such that $v_0+v_1-\sum_2^rs_iv_i=0$, given by the fact that the toric variety is non singular.

Let us now list the invariant rational curves on $\ex=\Proj{\el{0}\oplus...\oplus\el{k}}$. They correspond to $(r-1)$-dimensional cones in the fan $\Delta$, as described in \ref{subvarieties}:
\begin{itemize}
\item\underline{pull-backs $V(\rho^i)$.} Every $(m-1)$-dimensional cone $\rho\in\Delta'$ determines $k$ $(r-1)$-dimensional cones in $\Delta$. If $\rho=\oplus_2^m\real_+\cdot n_j$,
$$\rho^i=\left(\bigoplus_{j=2}^m\real_+\cdot p(n_j)\right)\bigoplus\left(\bigoplus_{l\neq i}\real_+\cdot\e_l\right)$$
\item\underline{fiber-curves $V(\rho^{\sigma}_{ij})$. } For every $i\neq j$ and every fixed point $V(\sigma)$
on $\exy$,  where $\sigma=\oplus_l\real_+\cdot n_l\in\Delta'(m)$ there is an $(r-1)$-dimensional cone in $\Delta$:
$$\rho_{ij}=\left(\bigoplus_l(\real_+\cdot p(n_l)\right)\bigoplus(\left(\bigoplus_{t\neq i,j}\real_+\cdot\e_t\right)$$
\end{itemize}

Consider an $(r-1)$-dimensional cone of the form $V(\rho^i)$. Because $\ex$ and $\exy$ are non singular there are integers $s_i,\s_i$ and $h_j$ such that:
\begin{equation}\label{s}
n_0+n_1-s_2n_2-...-s_m n_m=0
\end{equation}
\begin{equation}\label{sbar}
p(n_0)+p(n_1)-\s_2n_2-...-\s_m n_m-\left(\sum_{i\neq j}h_i\e_i\right)=0
\end{equation}
Using the definition of the map $p$ we get:
$$n_0+n_1-\s_2n_2-...-\s_m n_m+\e_0(a^0_0+a^0_1-\sum\s_j a^0_j)+...+\e_k(a^k_0+a^k_1-\sum\s_j a^k_j)-(+_{i\neq j}h_i\e_i)=0$$
The relation $\e_0+...+\e_k=0$ and $a^i_0+a^i_1-\sum\s_j a^i_j=\el{i}\cdot\rho$ yields:
$$ n_0+n_1-\s_2n_2-...-\s_m n_m-\sum_{j\neq i}(\e_j(\el{j}\cdot\rho-\el{i}\cdot\rho-h_j))=0$$ from which we conclude that
\begin{itemize}
\item $s_i=\s_i$ for $i=2,...,m$;
\item $h_j=(\el{j}-\el{i})\cdot\rho$.
\end{itemize}
It follows that:
\vspace{.1in}

\begin{tabular}{|c|}
\hline
$V(\e_0)\cdot V(\rho^i)=0$ if $i=0$\\
$V(\e_0)\cdot V(\rho^i)=(\el{i}-\el{0})\cdot\rho$ if $i\neq 0$\\
 $\pi^*(H)\cdot V(\rho^i)=H\cdot \rho$ for every line bundle $H$ on $\exy$.\\
\hline
\end{tabular}
\vspace{.1in}

Let us now consider a fiber-curve $V(\rho_{ij}^{\sigma})$ and integers $s_i, h_i$ such that
$$\e_i+\e_j-\left(\sum_{l\neq i,j}h_l\e_l\right)-\left(\sum_1^m s_t p(n_t)\right)=0$$
Because the $n_i$'s can be chosen to form a $\zed$-basis for $\Delta'$ it is clear that
$$ h_l=-1\;\text{ for }l\neq i,j\;\;\;\text{ and }s_t=0\;\;\text{ for }t=1,...,m$$
which implies that
\vspace{.1in}

\begin{tabular}{|c|}\hline
 $V(\e_0)\cdot V(\rho_{ij}^{\sigma})=1$ for all $i,j$;\\
$\pi^*(H)\cdot V(\rho_{ij}^{\sigma})=0$ for every line bundle $H$ on $\exy$.\\
\hline
\end{tabular}
\vspace{.1in}
\subsection{The polarization}
We will now describe the polarization associated to a defect toric manifold, i.e. the line bundle $\oof{\ex}{1}$ defining the embedding in projective space. From the description of the Picard group we can write
$$\oof{\ex}{1}=a\xi+\pi^*(H)$$ where $H$ is a line bundle on the non singular toric variety $\exy$ and $a$ is an integer.
In \ref{dual} we have seen that $\oof{\ex}{1}|_F\cong\oof{\proj{r-m}}{1}$ for a general fiber of the morphism $\pi:\ex\to\exy$. This implies that for every fiber-curve $V(\rho_{ij})$
$$\oof{\ex}{1}\cdot V(\rho_{ij}^{\sigma})=a(\xi\cdot V(\rho_{ij}^{\sigma}))+\pi^*(H)\cdot V(\rho_{ij}^{\sigma})=$$
$$=a(V(\e_0)+\pi^*(\el{0}+H))\cdot V(\rho_{ij}^{\sigma})=a=1$$
We can now see that:
\begin{proposition}\label{projective}
Let $\ex=\Proj{\el{0}\oplus...\oplus\el{k}}$ be a toric defect manifold polarized by the line bundle $\oof{\ex}{1}=\xi\otimes\pi^*(H)$. Then the line bundles $\el{i}\otimes H$ are very ample on $\exy$.
\end{proposition}
\begin{proof}Recall that a line bundle $\el{}$ on a non singular toric variety is very ample if and only if $\el{}\cdot V(\rho)=h$, where $\el{}|_{V(\rho)}=\oof{\proj{1}}{h}$, with $h\geq 1$ for  all the invariant rational curves $V(\rho)$, see for example \cite{dr}.
Using the intersection numbers estimated in the previous subsection we see that:
\begin{itemize}
\item $\oof{\ex}{1}\cdot V(\rho^i)=(V(\e_0)+\pi^*(\el{0})\otimes H)\cdot V(\rho^i)=(\el{0}\otimes H)\cdot \rho\geq 1$ if $i=0$;
\item $\oof{\ex}{1}\cdot V(\rho^i)=(V(\e_0)+\pi^*(\el{0})\otimes H)\cdot V(\rho^i)=(\el{i}\otimes H)\cdot \rho\geq 1$ if $i\neq 0$;
\item $\oof{\ex}{1}\cdot V(\rho_{ij}^{\sigma})=(V(\e_0)+\pi^*(\el{0})\otimes H)\cdot V(\rho_{ij}^{\sigma})=1$ for every $ij$.
\end{itemize}
It follows that for every invariant rational curve $\rho$ in $\Delta'$,  $(\el{i}\otimes H)|_{V(\rho)}=\oof{\ex}{1}|_{V(\rho^i)}$ for $i=0,...,k$. Hence the line bundles $\el{i}\otimes H$ are very ample on the $\exy$. \end{proof}


\section{Defect Polytopes}\label{defpoly}

In this section we will use geometry to investigate combinatorial properties of polytopes.

Let $N$ be a $d$-dimensional lattice and $V=N\otimes \real$. Throughout this paper the word polytope will always mean an {\em integral convex simple polytope} $P\subset V$, i.e. vertices of $P$ lie in $N$ and any vertex is incedent to exactly $\dim(P)$ edges. We will denote by $P(k)$ the set of $k$-dimensional faces of $P$.

Let $\dim(P)=r$ and let ${\mathcal A}_P$ denote the affine span of $P$. Then there is an invertible affine transformation $\phi:{\mathcal A}_P\to \real^d$ with $\phi({\mathcal A}_P\cap V)=\zed^r$. The volume of $P$ is defined as 
$${\mathrm Vol}(P)={\mathrm Leb}(\phi(P))$$
where ${\mathrm Leb}$ denotes the Lebesgue measure on $\real^r$. Similarly for each face $F\subset P$. In particular 
\begin{itemize} 
\item For every edge $\ell\in P(1)$ the length $\mathrm{Vol}(\ell)=|\ell\cap N|-1$. 
\item Let $\Delta_r$ denote the $r$-dimensional standard simplex. Then  $|\Delta_r(k)|={r+1\choose k+1}$  and $Vol(F)=\frac{1}{k!}$ for each face $F\in \Delta_r(k)$.
\end{itemize}
Consider the combinatorial invariant defined in (\ref{c(P)}).
$$c(P)=\sum_{F\subset P}(-1)^{\mathrm {codim}(F)}(\dim(F)+1)! \mathrm{Vol}(F)$$
The equality 
$$\sum_{k=0}^{r}(-1)^{r-k}(k+1)!{r+1\choose k+1}\frac{1}{k!}=0$$
yields that 
$$c(\Delta_r):=\sum_{k=0}^{r}(-1)^{r-k}(k+1)!\left(\sum_{F\in \Delta_r(k)}Vol(F)\right)=0$$

Recall that there is a one to one correspondence between polytopes and projective toric varieties. We will denote by $X(P)$ the projective toric variety associated to $P$, and by $\oof{P}{1}$ the invertible sheaf defining the embedding in $\Pin{|P\cap N|-1}$. The complete fan defining the toric variety $X(P)$ will be denoted by $\Delta(P)$. Recall that the variety $X(P)$ is Cohen-Macaulay, \cite[3.9]{Oda}. Geometry and combinatorics are related by the following correspondence:
\begin{center}
\begin{tabular}{lllll}

 $F_{\sigma}\in P(k)$&$\longleftrightarrow $&$V(F_{\sigma})\subset X(P)$ &$\longleftrightarrow$ &$\sigma_F\in\Delta_P(r-k)$\\
&&invariant of dimension $k$&&
\end{tabular}
\end{center}

Let $0\in V(F_{\sigma})$ be the unique $0$-dimensional toric orbit. Then 
$$\mult_F= \mult_0( V(F_{\sigma}))$$ measures how singular is $X(P)$
along $ V(F_{\sigma})$. Equivalently \cite[8.2]{D} if
$\sigma=<e_1,...,e_k>$ then $\mult_F$ is equal to the number of lattice points in the parallelotope $\{\sum_i\alpha_i e_i\,|\, 0\leq\alpha_i<1\}$.

An $r$-dimensional polytope is {\em Delzant} if every vertex $m$ is incident to exactly $r$ edges $\{l_1,...,l_r\}$ and $\{m_1-m,...,m_r-m\}$ form a $\zed$-basis of $N$, where $m_i$ are the lattice points on $l_i$ next to $m$. 
 The polytope $P$ is Delzant if and only if $X(P)$ is a non singular toric variety, i.e. a toric manifold, \cite[2.22]{Oda}. Clearely $X(P)$ is non singular if and only if $\mult_F=1$ for every $F\in P(k)$, $k=0,...,r$.
\subsection{Jets on $X(P)$} 

We shall recall briefly the definition and basic properties of the first jet bundle $J_{1}(\oof{P}{1})$, or sheaf of principal parts, associated to  $\oof{P}{1}$. For details we refer to \cite{LaTh, book, dr}.

Consider the projections $\pi_i:X(P)\times X(P)\to X(P)$, $i=1,2$,  and the ideal sheaf ${\mathcal I}_{\Delta}$ defining the diagonal $\Delta$. This data defines a sheaf of principal parts, supported on the reduced subscheme  $\Delta$.
$$J_{1}(\oof{P}{1})=\pi_{2_*}(\pi_1^*(\oof{P}{1}\otimes\ool{X(P)\times X(P)}{}/{\mathcal I}_{\Delta}^2)$$ 
The following  reports some  properties classically known, \cite[\S4]{LaTh}.

Let $i:M\to X(P)$ be the inclusion of the non singular scheme and $\Omega^1_M$ the sheaf of $1$-differential forms. The sheaf of $1$-differential form on $X(P)$ is $\Omega^1_P=i^*\Omega^1_M$.
\begin{lemma}\label{jetseq}
There is an exact sequence called the {\it first jet sequence}:
$$
0\to \Omega^1_P\otimes\oof{P}{1}\to J_1(\oof{P}{1})\to \oof{P}{1}\to 0 $$ 
\end{lemma}

\begin{proof}
Consider the exact sequence
$$0\to {\mathcal I}_{\Delta}/{\mathcal I}_{\Delta}^2\to {\mathcal O}/
{\mathcal I}_{\Delta}^2\to {\mathcal O}/{\mathcal I}_{\Delta}\to 0$$
tensoring by $\pi_1^*(\oof{P}{1})$ and pushing forward by $\pi_2$ gives the result, considering that \\
$R^1\pi_*({\mathcal I}_{\Delta}/{\mathcal I}_{\Delta}^2)=0$ and $\pi_{2_*}(\pi_1^*(\oof{P}{1})\tensor {\mathcal I}_{\Delta}/{\mathcal I}_{\Delta}^2)\cong \Omega_P^1(1)$.
\end{proof}
Let $r=\dim(P)=\dim(X(P))$. A corollary of the previous Lemma is that 
$J_1(\oof{P}{1})$ is locally free of rank $(r+1)$.
The first jet of a section $s\in H^0(X,\oof{P}{1})$ at a smooth point $x$ is $j_1(x,s)=(a_0,...,a_r)$, defined by the coefficients of the Taylor expansion of $s$ about $x$, truncated to the first degree.

Because $H^1(X(P), \Omega^1_X\otimes\oof{P}{1})=0$  Lemma \ref{jetseq} implies that the global sections $H^0(X(P),\oof{P}{1})$ can be naturally identified with a subvector space of $H^0(X(P), J_1(\oof{P}{1})$.

Hence there is a natural map of vector bundles:
\begin{equation}
j_1:X\times H^0(X, \oof{P}{1})\to J_1(\oof{P}{1}) \label{jetmap}
\end{equation}

On smooth points the map  assigns to each section its first jet. On a closed point $x$ the map corresponds to
$$ H^0(X, \oof{P}{1})\to H^0(x,\oof{P}{1}\otimes{\mathcal O}_X/m_x^2)$$

When $P$ is Delzant, i.e. X(P) is non singular, the kernel of this map is the conormal bundle shifted by one, $N^*_{X(P)/\proj{|P\cap N|-1}}(1)$, and the image of the induced map 
$${\mathbb P}(N_{X(P)/\proj{|P\cap N|-1}}(-1))\to X(P)\times({\proj{|P\cap N|-1}})^*\to ({\proj{|P\cap N|-1}})^*$$ defines the dual variety $\exd$. The dual defect is in fact controlled by the vanishing of Chern classes of the first jet bundle:
\begin{lemma}\label{vanishing}\cite[1.6.10]{book}
Notation as in section \ref{dual}. An $r$-dimensional toric manifold has dual defect $d>0$ if and only if $c_k(J_1(\oof{P}{1})=0$ for $k\geq r-d+1$ and $c_{r-d}(J_1(\oof{P}{1})\neq 0$.
\end{lemma}
In particular $d>0$ if and only if $c_r(J_1(\oof{P}{1})=0$.

The following proposition gives the key link between combinatorics and geometry, that will allow us to understand the behavior of the integer $c(P)$, based on the dual behavior of the variety $X(P)$.
Let
$$c^*(P)=\sum_{k=0}^{r}(-1)^{r-k}(k+1)!\left(\sum_{F\in
    P(k)}\frac{Vol(F)}{\mult_{F_{\sigma}}}\right)$$
Notice that in the non singular case $c^*(P)=c(P)$. The following
    Proposition showes that the combinatorial invariant $c^*(P)$ is an integer.
\begin{proposition}\label{main}
Let $P$ be a  polytope. Let $c_{r}(J_{1}(\oof{P}{1})$ be the top Chern class of the vector bundle $J_{1}(\oof{P}{1})$. Then
$$c^*(P)=c_{r}(J_{1}(\oof{P}{1}))$$

\begin{proof}
From (\ref{jetseq}) one sees that  
$$c_r(J_{1}(\oof{P}{1}))=\sum_{k=0}^r (r+1-k)c_k(\Omega^1_P)\cdot\oof{P}{1}^{r-k}$$
Consider the generalized Euler sequence, \cite[12.1]{BaCo}:
$$0\to \Omega^1_P\to\bigoplus_{\xi\in \Delta_P(1)} {\mathcal O}_{X(P)}(-V(\xi))\to{\mathcal O}^{|M\cap P(r-1)|-r}\to 0$$ 
The Chow groups $A_k(X(P))$ are generated by the classes $[V(\sigma)]$, for $\sigma\in\Delta(r-k)$, \cite[5.1]{Fu}.
The total Chern class of $\Omega^1_P$ is then
$$c(\Omega^1_P)=\prod_{\xi\in \Delta_P(1)}(1-[V(\xi)])\,\,\,\,\text{ and }$$

$$c_k(\Omega^1_P)=(-1)^k\sum_{\xi_1\neq...\neq \xi_k\in\Delta(1)} [V(\xi_1)]\cdot ...\cdot [V(\xi_k)]$$
Using the fact that: \cite[10.9]{D}
$$[V(\xi_1)]\cdot ...\cdot [V(\xi_k)]=\frac{1}{\mult_{F_{\sigma}}}[V(F_{\sigma})],\,\, \sigma=<\xi_1,...,\xi_k>$$
$$c_k(\Omega^1_P)=(-1)^k\sum_{\sigma\in\Delta_P(k)}\frac{1}{\mult_{F_{\sigma}}}[V(F_{\sigma})]=\sum_{F_{\sigma}\in P(r-k)}(-1)^k\frac{1}{\mult_{F_{\sigma}}}[V(F_{\sigma})]$$
where $F_{\sigma}\in P(k)$ is the $k$-dimensional face dual to the cone $\sigma_F\in\Delta_P(r-k)$. 

Moreover \cite[5.3]{Fu}:
\begin{equation}\label{integer}
\mathrm{Vol}(F_{\sigma})=\frac{{\mathrm deg}(\oof{P}{1}^k\cap [V(F_{\sigma})])}{k!}
\end{equation}

for example $\mathrm{Vol}(P)=\frac{deg(\oof{P}{1})}{r!}$ and for any $r$-dimensional cone $\sigma\in\Delta(r)$, $V(F_{\sigma})$ is a fixed point of the toric action, and $F_{\sigma}=m(\sigma)$ is a vertex of the polytope $P$. So in this case
$$\frac{{\mathrm deg}(\oof{P}{1}^k\cap [V(F_{\sigma})])}{k!}={\mathrm deg}(V(f_{\sigma}))=1$$
 It follows that
$$c_r(J_{1}(\oof{P}{1}))=\sum_{k=0}^r (r+1-k)(-1)^k\left(\sum_{F_{\sigma}\in
P(r-k)}\frac{1}{\mult_{F_{\sigma}}}[V(F_{\sigma})]\right)\cdot\oof{P}{1}^{r-k}$$
and thus
$$c_r(J_{1}(\oof{P}{1}))=\sum_{k=0}^r (r+1-k)(-1)^k(r-k)!\left(\sum_{F_{\sigma}\in
P(r-k)}\frac{Vol(F_{\sigma})}{\mult_{F_{\sigma}}}\right)=$$   
$$=\sum_{k=0}^{r}(-1)^{r-k}(k+1)!\left(\sum_{F\in P(k)}\frac{Vol(F_{\sigma})}{\mult_{F_{\sigma}}}\right):=c^*(P)$$
\end{proof}
\end{proposition}
The non-negativity of the integer $c(P)$ is a simple consequence of this identification.
\begin{corollary}\label{nonneg}
If $P$ is a (integral convex simple) polytope then $c^*(P)$ is a
nonnegative integer.
\end{corollary}
\begin{proof}
 The fact that the line bundle $\oof{P}{1}$ defines an embedding implies that the map $$ H^0(X, \oof{P}{1})\to H^0(x,\oof{P}{1}\otimes{\mathcal O}_X/m_x^2)$$ is surjective for every closed point $x$. It follows that the vector bundle map (\ref{jetmap}) is surjective on stalks and hence it is a surjection of vector bundles. This means that 
the first jet bundle is generated by $ H^0(X(P),\oof{P}{1})$, a subspace of its global sections. Any globally generated vector bundle has non negative top Chern class, \cite[12.1.7]{Fu2}. This proves that $c^*(P)\geq 0$. 
\end{proof}
\subsection{The classification of defect polytopes}

In this section we will classify Delzant polytopes ``minimal'' with respect to $c(P)$.
\begin{definition}\label{def1}
A Delzant polytope is said to be a {\em defect Delzant polytope} if $c(P)=0$.
\end{definition}

Obviously in dimension one there is no defect Delzant polytope.

In dimension $2$ there are no defect polytopes rather than the
$2$-dimensional simplex. This is quickly checked using Pick's formula:
\begin{equation}
|P\cap N|=\mathrm {Area}(P)+\frac{1}{2}\cdot \mathrm{Perimeter}(P)+1 \label{pick}
\end{equation} 
 If $r=2$ then
$$c_2(P)=6\cdot \mathrm{Area}(P)-2\cdot \mathrm{Perimeter}(P)+|P(0)|=$$
$$=6|P|-5\cdot \mathrm{Perimeter}+|P(0)|-6\geq 0$$
Because  $|P\cap N|\geq \mathrm{Perimeter}(P)$ and $|P\cap N|\geq 4, |P(0)|\geq 4$ unless $P$ is a
simplex in which case $|P\cap N|=|P(0)|=3$ and thus $c(P)=0$.

In higher dimension this is no longer the case. See for example the polytope in Figure \ref{poly}.
\begin{definition} \label{join}Let $P_0,...,P_k$ be $m$-dimensional polytopes having the same combinatorial type. Let $P_i={\mathrm Conv}\{m_1^i,...,m_s^i\}$ and let $(e_1,...,e_k)$ be a basis for $\zed^k$ and $e_0=0$. The {\em projective join} of $P_0,...,P_k$ is the $(m+k)$-dimensional polytope defined as:
$${\mathbb P}(P_0,...,P_k)=\mathrm{Conv}\{(m^i_j,e_i)\}_{j=1,...,s, i=0,...,k}$$
\end{definition}

\begin{figure}
\includegraphics{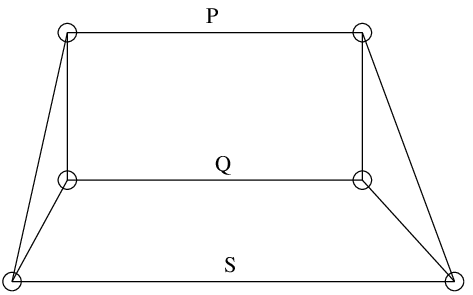}
\begin{caption}{\label{joinf} ${\mathbb P}(P,Q,S)$}
\end{caption}
\end{figure}
\begin{itemize}
 \item Figure \ref{joinf} shows the projective join of the three one dimensional
polytopes $P,Q,R$ of length $2,2,3$ respectively. 
\item Clearly a standard $r$-dimensional simplex is given by ${\mathbb P}(m,...,m)$ for any $0$-dimensional polytope $m$.
\item The two dimensional projective joins are the standard simplex and the $4$-polygons with two parallel edges and two edges of length $1$.
\item The product of a $k$-dimensional standard simplex and an $m$-dimensional polytope $Q$ is given by the projective join ${\mathbb P}(Q,...,Q)$.      
\end{itemize}
\begin{proposition}\label{dualpoly}
Let $P$ be an $r$-dimensional Delzant polytope. $P$ is a defect polytope if and only if there exists a positive integer $k$, $\max(2,\frac{r+1}{2})\leq k\leq r$, and $(r-k)$-dimensional Delzant polytopes, $P_0,...,P_{k}$,  such that $P={\mathbb P}(P_0,...,P_{k})$.
\end{proposition}
\begin{proof} $P$ is a defect polytope if and only if $c(P)=0$. Because of \ref{vanishing} and \ref{main} $P$ is a defect polytope if and only if the toric manifold $X(P)$ has positive dual defect. From \ref{dual} it follows that $X(P)=\Proj{\el{0}\oplus...\oplus\el{\frac{r+d}{2}}}$, where the $\el{i}$ are line bundles on a non singular toric variety $Y$, defined by the fan $\Delta'$ of dimension $(\frac{r-d}{2})$. The polytope $P$ is defined by the line bundle $\oof{P}{1}=\xi\otimes\pi^*(H)$. By \ref{projective} the line bundles $\el{i}\otimes H$ are very ample and hence define Delzant polytopes $P_0,...,P_{\frac{r+d}{2}}$. The polytope $P$ is constructed dually to the fan of $\Proj{\el{0}\oplus...\oplus\el{\frac{r+d}{2}}}$, described in the previous chapter. In particular the length of the edges are given by the intersection numbers of the line bundle $\oof{P}{1}$ with the invariant curves, as described in \ref{projective}. Let $P_i=\mathrm{Conv}\{m_1^i,...,m_s^i\}$, where!
 $s=|\Delta'(1)|=|P_i(\frac{r-d}{2
}-1)|$.
\begin{center} 
\begin{tabular}{lll}
$k$&$F\subset P(k)$& $\sigma_F\in\Delta_P(r-k)$\\
\hline
$0$&$m_{\sigma}^i$&  $p(\sigma)+\sigma^i,\sigma\in\Delta'(\frac{r-d}{2})$\\
$1$&edge connecting $m_{\sigma_1}^i$ and $m_{\sigma_2}^i$&$\rho^i$, $\rho=\sigma_1\cap\sigma_2$\\
&edge connecting $m_{\sigma}^i$ and $m_{\sigma}^j$&$\rho_{ij}^{\sigma}$\\
$\frac{r-d}{2}$& $P_i$&$p(0)+e^i$\\
$\frac{r+d}{2}$& $\Delta_{\frac{r+d}{2}}$&$p(\sigma)$\\

\end{tabular}
\end{center}
\vspace{.1in}

It is straight forward to check that 
$$P={\mathbb P}(P_0,...,P_{\frac{r+d}{2}})=\mathrm{Conv}\{m_{\sigma}^i\}_{\sigma\in\Delta'(\frac{r-d}{2}), k=0,...,\frac{r+d}{2}}$$
Set $k=\frac{r+2}{2}$. Because $r+d\geq 4$ and  $d\geq 0$ we have $\max(2,\frac{r+1}{2})\leq k\leq r$.

 For example the polytope in Figure \ref{joinf} defines the toric manifold ${\mathbb P}(\oof{\proj{1}}{1}\oplus\oof{\proj{1}}{1}\oplus\oof{\proj{1}}{2})$.
\end{proof}
It follows that:
\begin{itemize}
\item The only defect polytope in dimension $2$ is the standard simplex ${\mathbb P}(m,m,m)$, as observed earlier.
\item The only defect Delzant polytopes in dimension three are the standard simplex ${\mathbb P}(m,m,m,m)$ and  ${\mathbb P}(P_0,P_1,P_2)$, where $P_0,P_1,P_2$ are $1$-dimensional polytopes.
\item the only Delzant defect four dimensional polytopes are the standard simplex  ${\mathbb P}(m,m,m,m,m)$ and ${\mathbb P}(P_0,P_1,P_2, P_3)$, where $P_0,P_1,P_2, P_3$ are $1$-dimensional Delzant.
\end{itemize}
One could give an explicit classification for any dimension.
\section{Final remarks and conjectures}\label{conjsec}
\subsection{The degree of the dual variety}

If $X(P)$ is non singular and has no dual defect, then the dual
variety is an hypersurface defined by a polynomial of degree equal to
the top chern class of the first jet bundle. Hence 
$$\mathrm{deg}(X(P)^*)=c(P)$$

This is the way this integer was defined in \cite{GKZ}.
More generally in \cite{GKZ} it is  shown that it corresponds to the degree of a rational homogeneous function, $D_P$,  and that it is nonnegative for any simple polytope, \cite[11,1.6]{GKZ}.

Several numerical experiments suggest that the invariant $c(P)$ should be nonnegative for any polytope. The following example shows that even when the function $D_P$ is not polynomial, the invariant $c(P)$ is nonnegative. 
\begin{examp}
Let $e_1,...,e_6$ be a $\zed$-basis for $\real^6$. Consider the $5$-dimensional hypersimplex 
$$\Delta(3,6)=\mathrm{Conv}\{(e_i+e_j+e_k), 1\leq i<j<k\leq k\}$$

The rational function $D_{\Delta(3,6)}$ is not polynomial, \cite[11, 2.5]{GKZ}. Let us compute the invariant $c(\Delta(3,6))$.
The number of faces are decribed in \cite[2.5]{FuSt}: 
\begin{center} 
\begin{tabular}{lllllll}
$k$&$0$&$1$&$2$&$3$&$4$&$5$\\
\hline
$|\Delta(3,6)(k)|$&$20$&$90$&$120$&$60$&$12$&$1$

\end{tabular}
\end{center}
\vspace{.1in}
Every face is a hypesymplex given by two disjoint partitions of the set $\{1,...,n\}$. Using the formula, \cite{FoSt}
$$\mathrm{Vol}(\Delta(k,n))=\frac{A_{n,k}}{(n-1)!}$$
where the numbers $A_{n,k}$ are the Eulerian numbers, we get
$$c(\Delta(3,6))= 352$$
\end{examp}

We conjecture that:
\begin{con}
the invariant $c(P)$ is nonnegative for any integral polytope.
\end{con}

\subsection{Other invariants}
The invariant $c(P)$ may seem not the most natural object to consider.
Define:
$$c_t(P)=\sum_{F\subset P}(-1)^{\mathrm{codim}(F)}(\mathrm{dim}(F)+t)!Vol(F)$$
The sum runs over all nonempty faces, including the polytope $P$.
The invariant $c_0(P)$ is certainly combinatorially more attractive, but unfortunately negative in several examples.
Even for standard simpleces $\Delta_r$ one can easily see that 
\begin{itemize}
\item $ c_0(\Delta_r)= (-1)^r$
\item $ c_i(\Delta_r)= 0$ for $i=1,...,r$;
\item $ c_i(\Delta_r)>0$ for $i\geq r+1$
\end{itemize}

The invariant $c(P)=c_1(P)$ is in fact the first for which one could have hoped to prove a non-negativity result.
Moreover it is a corollary of \ref{nonneg} that every $c_t(P)$ is nonnegative for $t\geq 1$.
\begin{corollary}
The invariant $c_t(P)\geq 0$ for $t\geq 1$, for any simple polytope $P$.
\end{corollary}
\begin{proof}
Proceed by induction on $t$. The assertion is true for $t=1$, by
\ref{nonneg} and \cite{GKZ}. Consider the polytope $P'=P\times [0,m]$.
$$c_t(P')=\sum _{F\subset P}[(-1)^{\mathrm{codim}(F)-1}(\mathrm{dim}(F)+t)!Vol(F)]+m[(-1)^{\mathrm{codim}(F)}(\mathrm{dim}(F)+t+1)!\mathrm{Vol}(F)$$

In the limit, as $m$ goes to infinity, the predominant term is $m c_{t+1}(P)$. It follows that $c_{t+1}(P)$ is non negative if $c_t(P)$ is non negative.\end{proof}

If $F$ is a face of an $r$-dimensional polytope $P$, let $\zed^F=(\text{Affine span of }F)\cap \zed ^r.$ Consider the following function of $n$:
$$f(P,n)=\sum_0^r (-n)^{r-k}(k+1)!\sum_{F\in P(k)}|\zed^F\cap nF|=\sum_0^r d_i(P) n^i$$
If $\dim(F)=k$ then the coefficient of the $n^k$ term in $|\zed^F\cap nF|$ is $\mathrm{Vol}(F)$. This implies that 
$$ c(P)=d_r(P)$$

Because of (\ref{integer}) $f(P,n)$ is a polynomial in $n$ with integer coefficients.
Some numerical experiments show evidence that every coefficient of this polynomial is nonnegative.
\begin{con}
$f(P,n)$ is a polynomial in $n$ with nonnegative integer coefficients.
\end{con}
Equivalently, we conjecture that for a  projective variety $X(P)$
$$f(P,n)=\sum_0^r (-n)^{r-k}(k+1)!\sum_{F\in P(k)}dim(H^0(\oof{F}{n}))$$ is a a polynomial in $n$ with nonnegative integer coefficients.
Because the line bundle $\oof{F}{1}$ is very ample for each face $F$ this would mean that
\begin{equation}\label{conj}
f(P,n)=\sum_0^r (-n)^{r-k}(k+1)!\sum_{F\in P(k)}\chi(\oof{F}{n}))
\end{equation}
is a polynomial in $n$ with nonnegative integer coefficients.

One could give one additional interpretation of (\ref{conj}). Let $C_n^m$ denote the $m$-dimensional cube with edges of length $n$. Geometrically it corresponds to the the Segre embedding  of $m$ copies of the $n$-th Veronese embedding of ${\mathbb P}^1:$
$$V_n^m:={\mathbb P}^1\times ...\times {\mathbb P}^1\to {\mathbb P}^{(n+1)^m-1}$$
  For any face $F$ consider the polytope $F\times  C_n^m$ associated to the variety $V(F)\times V_n^m$. Let $\pi_i$ be the projection onto the $i$-th factor, $i=1,2$. The polytope $F\times  C_n^m$ is defined by the invertible sheaf:
$$\pi_1^*(\oof{F}{1})\otimes\pi^*_2(\oof{V_n^m}{1})=\oof{V(F)\times V_n^m}{1}$$
We conjecture that 
$$f(P,n)=\sum_F (-1)^{\mathrm{codim}(F)}(\mathrm{dim}(F))!\chi(\oof{V(F)\times V_n^m}{1})$$
is a polynomial in $n$ with nonnegative integer coefficients.

For standard simpleces one can show the nonnegativity of the constant term because $$d_0(\Delta_r)=c_{r+1}(\Delta_r)$$
One can see this equality using the Chu-Vendermonde identity
\cite{RR}.
\subsection{Defect polytopes}

We were not able to find examples of polytopes with  $c(P)=0$ and which do not have the structure of a projective join.
\begin{con}
$P$ is a defect polytope if and only if it is a projetive join of some kind.
\end{con}
Our approach unfortunately does not work  in a more general setting. The use of duality and intersection estimates relies strongly on the fact that the variety is nonsingular.

A purely combinatorial proof of the results contained in this paper could possibly open the door to a generalization.



\begin{thebibliography}{10}
\bibitem[BaCox]{BaCo} 
Batyrev, Victor V.; Cox, David A.
\newblock {\it  On the Hodge structure of projective hypersurfaces in toric varieties}.  
\newblock Duke Math. J.  {\bf 75}  (1994),  no. 2, 293--338.
\bibitem[BeSo]{book}
Beltrametti, Mauro C.; Sommese, Andrew J. 
\newblock {\it The adjunction theory of complex projective varieties.} 
\newblock de Gruyter Expositions in Mathematics, {\bf 16}. Walter de Gruyter \& Co., Berlin, 1995.
\bibitem[BeFaSo]{bfs}
Beltrametti, Mauro C.; Fania, M. Lucia; Sommese, Andrew J. 
\newblock {\it On the discriminant variety of a projective
manifold.} 
\newblock Forum Math. {\bf 4} (1992), no. 6, 529--547.

\bibitem[DR]{dr}Di Rocco, Sandra.
\newblock {\it Generation of $k$-jets on toric varieties.}
\newblock  Math. Z.  {\bf 231}  (1999),  no. 1, 169--188
\bibitem[DRSo]{chern}  Di Rocco, Sandra, Sommese, Andrew, J.
\newblock {\it Chern numbers of ample  vector bundles on toric surfaces}.  
\newblock Transactions of the A.M.S., to appear (AG9911192 ).
\bibitem[E1]{E1}Ein, Lawrence.
\newblock {\it Varieties with small dual varieties. II}.  
\newblock Duke Math. J.  {\bf 52}  (1985),  no. 4, 895--907.
\bibitem[E2]{E2}Ein, Lawrence.
\newblock {\it Varieties with small dual varieties. I}. 
\newblock  Invent. Math.  {\bf 86}  (1986),  no. 1, 63--74.
\bibitem[FoSt]{FoSt}  Foata, Dominique 
\newblock {\it Distributions eule'riennes et mahoniennes sur le groupe des permutations.} (French) With a comment by Richard P. Stanley.
\newblock  NATO Adv. Study Inst. Ser., Ser. C: Math. Phys. Sci., {\bf 31},  Higher combinatorics (Proc. NATO Advanced Study Inst., Berlin, 1976),  pp. 27--49, Reidel, Dordrecht-Boston, Mass., 1977 
\bibitem[FU]{Fu} Fulton, William.
\newblock {\it Introduction to toric varieties}. 
\newblock Annals of Mathematics Studies, {\bf 131}. The William H. Roever Lectures in Geometry. Princeton University Press, Princeton, NJ, 1993.
\bibitem[FU2]{Fu2} Fulton, William 
\newblock {\it Intersection theory.} Second edition. 
\newblock  Springer-Verlag, Berlin, 1998.
\bibitem[FuSt]{FuSt}Fulton, William; Sturmfels, Bernd.
\newblock {\it Intersection theory on toric varieties.}
\newblock  Topology  {\bf 36}  (1997),  no. 2, 335--353.
\bibitem[GKZ]{GKZ}Gelfand, I. M.; Kapranov, M. M.; Zelevinsky, A. V. 
\newblock{\it Discriminants, resultants, and multidimensional determinants.}
\newblock  Mathematics: Theory \& Applications. Birkha\"user Boston, Inc., Boston, MA, 1994.
\bibitem[Dan]{D}  Danilov, V. I. 
\newblock{\it The geometry of toric varieties.} 
\newblock Russian Math. Surveys {\bf 33}  (1978), no. 2(200), 85--134, 247.  
\bibitem[LaTh]{LaTh}  Laksov, Dan; Thorup, Anders
\newblock {\it Weierstrass points on schemes}.
\newblock  J. Reine Angew. Math.  {\bf 460}  (1995), 127--164. 
\bibitem[LaSt86]{ls1}
Lanteri, Antonio; Struppa, Daniele 
\newblock {\it Projective manifolds whose topology is strongly reflected in their hyperplane sections}. 
\newblock Geom. Dedicata {\bf 21} (1986), no. 3, 357--374.
\bibitem[LaSt87]{ls}
Lanteri, Antonio; Struppa, Daniele 
\newblock {\it Projective $7$-folds with positive defect}. 
\newblock Compositio Math. {\bf 61} (1987),
no. 3, 329--337.
\bibitem[Lib]{lib}
Lieberman D.I.
\newblock {\it Holomorphic vector fields and rationality}.
\newblock Group actions and vector fields (Vancouver, B.C., 1981),
99-117, lecture Notes in Math., {\bf 956}, Springer, 1982.
\bibitem[ODA]{Oda} Oda, Tadao.
\newblock  {\it Convex bodies and algebraic geometry. An introduction to the theory of toric varieties}.
\newblock Ergebnisse der Mathematik und ihrer Grenzgebiete (3)  {\bf 15}. Springer-Verlag, Berlin, 1988.
\bibitem[RR]{RR}  Roy, Ranjan 
\newblock {\it Binomial identities and hypergeometric series.}
\newblock  Amer. Math. Monthly  {\bf 94}  (1987),  no. 1, 36--46. 
\end{thebibliography}
 \end{document}